\documentclass[reqno]{amsart}

\usepackage{amsfonts,amssymb,amsmath,latexsym,indentfirst}

\usepackage[notcite,notref]{showkeys}

\newtheorem{prop}{Proposition}[section]
\newtheorem{rema}{Remark}[section]

\newtheorem{lemm}{Lemma}[section]
\newtheorem{theo}{Theorem}[section]

\newcommand{\R}{\ensuremath{{\mathbb{R}} }}

\newcommand{\peq}{\hspace*{0.10in}}
\newcommand{\peqq}{\hspace*{0.05in}}

\newcommand{\fim}{\rightline{$\blacksquare$}}

\newcommand{\V}{\vspace*{0.25in}}

\begin{document}
\title[Generalized Boussinesq equation]{Global solutions below the energy space for the generalized Boussinesq equation}
\footnotetext{The first author is partially supported by CNPq-Brazil and FAPEMIG-Brazil.}
\date{}

%\author[]{}
%\address{}
%\email{}
%\thanks{}
%\begin{abstract}
%\end{abstract}

\maketitle

\begin{center}
 {\large \textbf{Luiz G. Farah}}\\ {\small ICEx, Universidade Federal de Minas Gerais \\
Av. Ant\^onio Carlos, 6627, Caixa Postal 702, 30123-970, Belo Horizonte-MG, Brazil.    \\
E-mail: lgfarah@gmail.com}\\
 \end{center}

\begin{center}
 {\large \textbf{Hongwei Wang}}\\ {\small 
Faculty of Science, Xi'an Jiaotong University\\
Xi'an 710049, P.R.China\\
and\\
Department of Mathematics, Xinxiang College\\
Xinxiang 453003, P.R.China.\\
E-mail: wang.hw@stu.xjtu.edu.cn}\\
 \end{center}

\begin{abstract}
We show that the Cauchy problem for the defocusing generalized Boussinesq
equation $u_{tt}-u_{xx}+u_{xxxx}-(|u|^{2k}u)_{xx}=0$, $k\geq1$, on the real line is globally well-posed in
$H^{s}(\R)$ for $s>1-({1}/{3k})$. To this end we use the Ò$I$-methodÓ, introduced by J. Colliander, M. Keel, G. Staffilani, H. Takaoka and T. Tao \cite{CKSTT4, CKSTT3}, to define a modification of the energy functional that is Òalmost conservedÓ in time. Our result extends the previous one obtained by Farah and Linares \cite{LGF} for the case $k=1$.
\end{abstract}

\section{Introduction}
We study the following initial value problem for a defocussing generalized Boussinesq
equation
\begin{eqnarray}\label{NLB}
\left\{
\begin{array}{l}
u_{tt}-u_{xx}+u_{xxxx}-(|u|^{2k}u)_{xx}=0, \peq k \geq 1, \peq  x\in \R,\, t>0,\\
u(x,0)=\phi(x) ; \peqq \partial_tu(x,0)=\psi_x(x).
\end{array} \right.
\end{eqnarray}

Equations of this type model a large rang of physical phenomena such as dispersive wave propagation, nonlinear strings and shape-memory alloys (see, for instance, Boussinesq \cite{BOU}, Zakharov \cite{Z} and Falk \textit{et al} \cite{FLS}).

Natural spaces to study the initial value problem above are the classical Sobolev spaces $H^s(\R)$, $s\in \R$, which are defined via the spacial Fourier transform 
\begin{equation*}
\hat{f}(\xi)\equiv\int_{\R}e^{-ix\xi}f(x)dx,
\end{equation*}
as the completion of the Schwarz class $\mathcal{S}(\R)$ with respect to the norm
\begin{equation*}
\|f\|_{H^{s}(\R)}=\|\langle\xi\rangle^s \widehat{f}\|_{L^{2}(\R)}
\end{equation*}
where $\langle \xi \rangle = 1+|\xi|$.

Given initial datas $(\phi, \psi)\in H^{s}(\R)\times H^{s-1}(\R)$ and a positive time $T>0$, we say that a function $u: \R \times [0,T] \rightarrow \R$ is a real solution of \eqref{NLB} if $u \in C([0,T]; H^s(\R))$ and $u$ satisfies the integral equation
\begin{equation}\label{INT} 
u(t)= V_c(t)\phi+V_s(t)\psi_x+\int_{0}^{t}V_s(t-t')(|u|^{2k}u)_{xx}(t')dt',
\end{equation}
where the two operators that constitute the free evolution are defined via Fourier transform by the formulas
\begin{eqnarray*}
\widehat{V_c(t)}\phi(\xi)&=&\frac{e^{it\sqrt{{\xi}^2+{\xi}^4}}+ e^{-it\sqrt{{\xi}^2+{\xi}^4}}}{2}\widehat{\phi}(\xi)\\
\widehat{V_s(t)}{\psi_x}(\xi)&=&\frac{e^{it\sqrt{{\xi}^2+{\xi}^4}}- e^{-it\sqrt{{\xi}^2+{\xi}^4}}}{2i\sqrt{{\xi}^2+{\xi}^4}}\widehat{\psi_x}(\xi).
\end{eqnarray*}

In the case that $T$ can be taken arbitrarily large, we shall say the solution is global-in-time. Here, we focus our attention in this case. 

Concerning the local well-posedness question, several results have been obtained in the last years for the generalized Boussinesq equation \eqref{NLB} ( see Bona and Sachs \cite{BS}, Tsutsumi and Matahashi \cite{TM}, Linares \cite{FL} and Farah \cite{LG3, LG4}). As far as we know, one has local well-posedness in $H^s(\R)$ for all $s>1/2-1/k$ \cite{LG3}. The same holds for the focusing case, that is,  equation \eqref{NLB} with positive sign in front of the nonlinearity. Note that this is exactly the same range obtained by Cazenave and Weissler\cite{CaW} for the nonlinear Schr\" odinger equation
\begin{eqnarray*}
iu_{t}+u_{xx}-(|u|^{2k}u)=0.
\end{eqnarray*}
We should point out that, up to now, there is no result addressing the ill-posedness question for the equation \eqref{NLB} with general $k$, so it is an interesting open problem.

Next we turn attention to the global-in-time well-posedness problem. It is well know that generalized Boussinesq equation enjoy the following conserved energy
\begin{equation}\label{EC}
E(u)(t)=\frac{1}{2}\|u(t)\|^2_{H^1}+ \frac{1}{2}\|(-\Delta)^{-\frac{1}{2}}\partial_tu(t)\|^2_{L^2}
+ \frac{1}{2k+2}\|u(t)\|^{2k+2}_{L^{2k+2}}.
\end{equation}

The local theory proved in \cite{FL} together with this conserved quantity immediately yield global-in-time well-posedness of (\ref{NLB}) for initial data $(\phi, \psi)\in H^1(\R)\times L^2(\R)$.  We should mention that the situation is very different in the focusing case: solutions may blow-up in finite time for arbitrary initial data $(\phi, \psi)\in H^1(\R)\times L^2(\R)$, see for instance Liu \cite{LIU2} and Angulo and Scialom \cite{AS}. We will not deal with this case in present work.

Our principal aim in the present work is to loosen the regularity requirements on the
initial data which ensure global-in-time solutions for the initial value problem (\ref{NLB}). This question have been already investigated by Farah and Linares \cite{LGF} in the particular case where $k=1$. Their approach were based on the $I$-method, invented by the $I$-team: Colliander, Keel, Staffilani, Takaoka and Tao. Although for the generalized Boussinesq equation (\ref{NLB}) scaling argument does not work and there is no conservation law at level $L^2$, we also successfully applied this method, in its first generation, obtaining global solutions in $H^s(\R)$ with $s<1$ for all $k\geq 1$.

We shall mention that there exists other refined versions of the $I$-method also introduced by the $I$-team in the context of nonlinear dispersive equations (see, for instance,  \cite{CKSTT3} and \cite{CKSTT5}). This approaches  have been applied for the Nonlinear Sch\"odinger equation and generalized KdV equation
$$\partial_t u+\partial_x^3u+\partial_x(u^{k+1})=0,$$
sometimes leading to sharp global results \cite{CKSTT3}. However, since the generalized Boussinesq equation (\ref{NLB}) has two derivatives in time, it is not clear whether this refined approachs can be use to improve our global result stated in Theorem \ref{T1} below.

The basic idea behind the $I$-method is the following: when $(\phi, \psi)\in H^{s}(\R)\times H^{s-1}(\R)$ with
$s<1$ in (\ref{NLB}), the norm $\|\psi\|^2_{H^1}$ could be infinity, and so the
conservation law (\ref{EC}) is meaningless. To overcome this
difficulty, we introduce a modified energy functional which is also defined for less regular functions.
Unfortunately, this new functional is not strictly conserved, but we
can show that it is \textit{almost} conserved in time. When one is
able to control its growth in time explicitly, this allows to iterate
a modified local existence theorem to continue the solution to any
time $T$. 

Now we state the main result of this paper.

\begin{theo}\label{T1}
The initial value problem (\ref{NLB}) is globally well-posed in $H^s(\R)$ for all \peqq $1-\dfrac{1}{3k}<s<1$. Moreover the solution satisfies
\begin{equation}\label{pb}
\sup_{t\in[0,T]}\left\{\|u(t)\|^2_{H^{s}}+ \|(-\Delta)^{-\frac{1}{2}}\partial_tu(t)\|^2_{H^{s-1}}\right\}\leq C(1+T)^{\frac{1-s}{6ks-6k+2}+}
\end{equation}
where the constant $C$ depends only on $s$, $\|\phi\|_{H^{s}}$ and $\|\psi\|_{H^{s-1}}$.
\end{theo}

The plan of this paper is as follows. In the next section we introduce some notation and preliminaries. 
Section 3 describes the modified energy functional. In Section 4, we prove the almost conservation law.  Section 5 contains the variant of local well-posedness result and the proof of the global result stated in Theorem \ref{T1}.

\section{Notations and preliminary results}\label{S3}

We use $c$ to denote various constants depending on $s$. Given any positive numbers $a$ and $b$, the notation $a \lesssim b$ means that there exists a positive
constant $c$ such that $a \leq cb$. Also, we denote $a \sim b$ when, $a \lesssim b$ and
$b \lesssim a$. We use $a+$ and $a-$ to denote $a+\varepsilon$ and $a-\varepsilon$, respectively, for
arbitrarily small $\varepsilon>0$.

We use $\|f\|_{L^p}$ to denote the $L^p(\R)$ norm and $L^q_tL^r_x$ to denote the mixed norm
\begin{equation*}
\|f\|_{L^q_tL^r_x}\equiv \left(\int \|f\|_{L^r_x}^q dt \right)^{1/q}
\end{equation*}
with the usual modifications when $q =\infty$.\\

We define the spacetime Fourier transform $u(t, x)$ by
\begin{equation*}
\widetilde{u}(\tau,\xi)\equiv\int_{\R}\int_{\R}e^{-i(x\xi+t\tau)}u(t,x)dtdx.
\end{equation*}
Note that the derivative $\partial_x$ is conjugated to multiplication by $i\xi$ by the Fourier transform.

We shall also define $D$ and $J$ to be, respectively, the Fourier multiplier with symbol $|\xi|$ and $\langle \xi \rangle = 1+|\xi|$. 
Thus, the Sobolev norms $H^s(\R)$ is also given by
\begin{equation*}
\|f\|_{H^s}=\|J^sf\|_{L^2_x}.
\end{equation*}

To describe our well-posedness results we define the $X_{s,b}(\R\times\R)$ spaces related to our problem (see also Fang and Grillakis \cite{FG} and \cite{LG4}). 
\begin{equation*}
\|F\|_{X_{s,b}(\R\times\R)}=\|\langle|\tau|-\gamma(\xi)\rangle^b\langle\xi\rangle^s \widetilde{F}\|_{L^{2}_{\xi,\tau}},
\end{equation*}
where $\gamma(\xi)\equiv\sqrt{{\xi}^2+{\xi}^4}$.

These kind spaces were used to systematically study nonlinear dispersive wave problems
by Bourgain \cite{B} and Kenig, Ponce and Vega \cite{KPV1, KPV2}. Klainerman and Machedon \cite{KM} also used similar ideas in their
study of the nonlinear wave equation. The spaces appeared earlier in the study of propagation of singularity 
in semilinear wave equation in the works \cite{RR}, \cite{Be83} of Rauch, Reed, and M. Beals.

For any interval $I$ we define the localized $X_{s,b}(\R \times I)$ spaces by 
\begin{equation*}
\|u\|_{X_{s,b}(\R \times I)}=\inf\left\{\|w\|_{X_{s,b}(\R \times\R)}:w(t)=u(t) \textrm{ on }  I\right\}.
\end{equation*}

We often abbreviate $\|u\|_{X_{s,b}}$ and $\|u\|_{X^{I}_{s,b}}$, respectively, for 
$\|u\|_{X_{s,b}(\R\times\R)}$ and $\|u\|_{X_{s,b}(\R \times I)}$.

We shall take advantage of the Strichartz estimate (see Ginibre, Tsutusumi and Velo \cite {GTV} for this inequality in the context of the Scr\"odinger equation. For the spaces $X_{s,b}$ defined above it follows by the argument employed by \cite{LGF}) 
\begin{equation}\label{L81}
\|u\|_{L^q_{t}L^p_{x}}\lesssim \|u\|_{X_{0,\frac{1}{2}+}}, \peq \textrm{where} \peq \dfrac{2}{q}=\dfrac{1}{2}-\dfrac{1}{p}.
\end{equation}

Taking $p=q$, we obtain the spacial case
\begin{equation}\label{L8}
\|u\|_{L^6_{x,t}}\lesssim \|u\|_{X_{0,\frac{1}{2}+}}
\end{equation}
which interpolate with the trivial estimate
\begin{equation}\label{TE}
\|u\|_{L^2_{x,t}}\lesssim \|u\|_{X_{0,0}}
\end{equation}
to give
\begin{equation}\label{L6}
\|u\|_{L^4_{x,t}}\lesssim \|u\|_{X_{0,\frac{1}{2}+}}.
\end{equation}

We also use 
\begin{equation*}
\|u\|_{L^{\infty}_tL^2_{x}}\lesssim \|u\|_{X_{0,\frac{1}{2}+}},
\end{equation*}
which together with Sobolev embedding gives
\begin{equation}\label{LI}
\|u\|_{L^{\infty}_{x,t}}\lesssim \|u\|_{X_{\frac{1}{2}+,\frac{1}{2}+}}.
\end{equation}

We also have the following refined Strichartz estimate in the case of differing frequencies (see  \cite{LGF}, Bourgain \cite{B1}).

% \begin{lemm}\label{L2.2}
% Let $\psi_1, \psi_2 \in X^{\delta}_{0,\frac{1}{2}+}$ be supported on spatial frequencies
% $|\xi|\sim N_1,N_2$, respectively. Then, if $N_1\lesssim N_2$, one has
% \begin{equation}\label{BOUR}
% \|\psi_1\psi_2\|_{L^2(\R\times[0,\delta])}\lesssim \frac{1}{N_2^{\frac{1}{2}}}\|\psi_1\|_{X^{\delta}_{0,\frac{1}{2}+}} \|\psi_2\|_{X^{\delta}_{0,\frac{1}{2}+}}.
% \end{equation}
% \end{lemm}

\begin{lemm}\label{L3}
Let $\psi_1, \psi_2 \in X_{0,\frac{1}{2}+}$ be supported on spatial frequencies $|\xi_i|\sim N_i$, $i=1,2$.
If $|\xi_1|\lesssim \min\left\{|\xi_1-\xi_2|,|\xi_1+\xi_2|\right\}$ for all
$\xi_i\in \textrm{supp}(\widehat{\psi}_i)$, $i=1,2$, then
\begin{equation}\label{GRU}
\|(D^{\frac{1}{2}}_x\psi_1)\psi_2\|_{L^2_{x,t}}\lesssim \|\psi_1\|_{X^{\delta}_{0,\frac{1}{2}+}}
 \|\psi_2\|_{X^{\delta}_{0,\frac{1}{2}+}}.
\end{equation}
\end{lemm}

Inequalities of this kind have been also obtained under the assumption $|\xi_2|\gg|\xi_1|$ for both Nonlinear Schr\"odinger equation and KdV equation (see Ozawa and Tsutsumi \cite{OT}, Gr\"unrock \cite{Gr05} and also \cite{B1}). Note that this relation implies the hypothesis of the above lemma.

\section{Modified energy functional}

In this section we brifly describe the $I$-method scheme. Given $s<1$ and a parameter $N\gg 1$, define the multiplier operator
\begin{equation*}
\widehat{I_Nf(\xi)}\equiv m_N(\xi)\widehat{f}(\xi),
\end{equation*}
where the multiplier $m_N(\xi)$ is smooth, radially symmetric, nondecreasing in $|\xi|$ and
\begin{eqnarray*}
m_N(\xi)=\left\{
\begin{array}{l l }
1&, \textrm{ if } |\xi|\leq N,\\
\left(\dfrac{N}{|\xi|}\right)^{1-s}&, \textrm{ if } |\xi|\geq 2N.
\end{array} \right.
\end{eqnarray*}

To simplify the notation, we omit the dependence of $N$ in $I_N$ and denote it only by $I$. Note that the operator $I$ is smooth of order $1-s$. Indeed, we have
\begin{equation}\label{smo}
\|u\|_{H^{s_0}}\leq c\|Iu\|_{H^{s_{0}+1-s}}\leq cN^{1-s}\|u\|_{H^{s_0}}.
\end{equation}

We can apply the operator $I$ in the equation \eqref{NLB}, obtaining the following modified equation
\begin{eqnarray}\label{MNLB}
\left\{
\begin{array}{l}
Iu_{tt}-Iu_{xx}+Iu_{xxxx}-I(|u|^{2k}u)_{xx}=0, \peq  x\in \R,\,t>0,\\
Iu(x,0)=I\phi(x); \quad \partial_tIu(x,0)=I\psi_x(x).
\end{array} \right.
\end{eqnarray}

Moreover, applying the operator $(-\Delta)^{-\frac{1}{2}}$ to the above equation \eqref{MNLB}, multiplying the result by $(-\Delta)^{-\frac{1}{2}}\partial_tIu$ and then integrating by parts with respect to $x$, we obtain
\begin{equation*}
\frac{1}{2}\frac{d}{dt}\left(\|Iu(t)\|^2_{H^1}+ \|(-\Delta)^{-\frac{1}{2}}\partial_tIu(t)\|^2_{L^2}\right)
+\langle I\left(|u|^{2k}u\right), \partial_tIu\rangle=0.
\end{equation*}

On the other hand,
\begin{equation*}
\frac{d}{dt}\|Iu(t)\|^{2k+2}_{L^{2k+2}}=(2k+2)\int_{\R}|Iu|^{2k}Iu\partial_tIu.
\end{equation*}

Therefore
\begin{equation}\label{ACL}
\frac{d}{dt}E(Iu)(t)=\langle |Iu|^{2k}Iu-I\left(|u|^{2k}u\right), \partial_tIu\rangle.
\end{equation}

By the Fundamental Theorem of Calculus, we have
\begin{equation}\label{FTC}
E(Iu)(\delta)-E(Iu)(0)=\int_0^{\delta}\frac{d}{dt}E(Iu)(t')dt'.
\end{equation}

Therefore, to control the growth of $E(Iu)(t)$ we need to understand how the quantity \eqref{ACL} varies in time.

\section{Almost conservation law}

In this section we will establish estimates showing that the quantity $E(Iu)(t)$ is \textit{almost} conserved in time.

\begin{prop}\label{p4.1}
Let $s>1/2$, $N\gg 1$ and $Iu$ be a solution of \eqref{MNLB} on $[0, \delta]$ in the sense of Theorem \ref{t3.2}. Then the following estimate holds
\begin{equation}\label{CC}
\left|E(Iu)(\delta)-E(Iu)(0) \right|\leqslant
CN^{-2+}\|Iu\|_{X_{1,\frac{1}{2}+}}^{2k+1}\|(-\triangle)^{-\frac{1}{2}}\partial_tIu\|_{X_{0,\frac{1}{2}+}^{\delta}}.
\end{equation}
\end{prop}

Before proceeding to the proof of the above proposition, we would like to make an interesting remark. The exponent $-2+$ on the right hand side 
of (\ref{CC}) is directly tied with the restriction $s > 1-({1}/{3k})$ in our main theorem. If one could replace the increment $N^{-2+}$ by $N^{-\alpha+}$ for some $\alpha>0$ the argument we give in Section \ref{GWP} would imply global well-posedness of (\ref{NLB}) for all $s > 1-({\alpha}/{6k})$. 
\V

\noindent\textbf{Proof. } Applying the the Parseval formula to identity \eqref{FTC} and using \eqref{ACL}, we have
\begin{equation*}
E(Iu)(\delta)-E(Iu)(0)=
\end{equation*}
\begin{equation*}
=\int_0^{\delta}\int_{\sum_{i=1}^{2k+2}\xi_i=0}\left(1-\frac{m(\xi_2+\xi_3+\cdots+
\xi_{2k+2})}{m(\xi_2) m(\xi_3)\cdots
m(\xi_{2k+2})}\right)\widehat{\partial_t
Iu(\xi_1)}\widehat{Iu(\xi_2)}\cdots\widehat{Iu(\xi_{2k+2})}.
\end{equation*}

Therefore, our aim is to obtain the following inequality
\begin{equation*}
{\bf Term}\leqslant
N^{-2+}\|I\phi_1\|_{X^{\delta}_{0,\frac{1}{2}+}}\prod\limits_{i=2}^{2k+2}\|I\phi_i\|_{X^{\delta}_{1,\frac{1}{2}+}},
\end{equation*}
where
\begin{equation*}
{\bf Term}\equiv \left|\int_0^{\delta}\int_{\sum_{i=1}^{2k+2}\xi_i=0}\left(1-\frac{m(\xi_2+\xi_3+\cdots+\xi_{2k+2})}{m(\xi_2)
m(\xi_3)\cdots m(\xi_{2k+2})}\right)\widehat{\partial_t
Iu(\xi_1)}\widehat{Iu(\xi_2)}\cdots\widehat{Iu(\xi_{2k+2})}\right|
\end{equation*}
and $\ast$ denotes integration over $\sum_{i=1}^{2k+2}\xi_i=0$.

We estimate $\mathbf{Term}$ as follows. Without loss of generality, we assume the Fourier transforms of
all these functions to be nonnegative. First, we bound the symbol in the parentheses pointwise
in absolute value, according to the relative sizes of the frequencies involved. After that, the
remaining integrals are estimated using Plancherel formula, H\"older's inequality and Lemma
\ref{L3}. To sum over the dyadic pieces at the end we need to have extra factors
$N_j^{0-}$, $j=1,\cdots,2k+2$, everywhere.

We decompose the frequencies $\xi_j$, $j=1,\cdots,2k+2$ into dyadic blocks $N_j$. By the symmetry
of the multiplier
\begin{equation}\label{MULT}
1-\frac{m(\xi_2+\xi_3+\cdots+\xi_{2k+2})}{m(\xi_2)m(\xi_3)\cdots m(\xi_{2k+2})}
\end{equation}
in $\xi_2, \xi_3, \cdots, \xi_{2k+2}$, we may assume for the
remainder of this proof that
\begin{equation*}
N_2\geqslant N_3\geqslant \cdots \geqslant N_{2k+2}.
\end{equation*}

Also note that $\sum_{i=1}^{2k+2}\xi_i=0$ implies $N_1\lesssim N_2$. We
now split the different frequency interactions into several cases,
according to the size of the parameter $N$ in comparison to the
$N_i$.\\

{\bf Case 1:} $N\gg N_2$.\\

 In this case, the symbol (2) is identically zero and the desired bound holds trivially.\\

{\bf Case 2:} $N_2\gtrsim N \gg N_3$.\\

Since
 $\Sigma_{i=1}^{2k+2}\xi_i=0$, we have here $N_1\thicksim N_2$. By the
 mean value theorem
 $$\left|\frac{m(\xi_2)-m(\xi_2+\xi_3+\cdots+\xi_{2k+2})}{m(\xi_2)}\right|\lesssim\frac{\left|
 \nabla m(\xi_2)\cdot(\xi_3+\xi_4+\cdots+\xi_{2k+2})\right|}{m(\xi_2)}\lesssim\frac{N_3}{N_2}.$$
 This pointwise bound together with Lemma \ref{L3} and Plancherel's theorem
 yield
\begin{eqnarray*}
 {\bf Term} &\lesssim& \frac{N_1N_3}{N_2}\|I\phi_1I\phi_3\|_{L^2(\mathbb{R}\times
 [0,\delta])}\|I\phi_2I\phi_4\|_{L^2(\mathbb{R}\times[0,\delta])}\prod\limits_{j=5}^{2k+2}\|I\phi_j\|_{L^{\infty}(\mathbb{R}\times[0,\delta])}\\
 &\lesssim&\frac{N_1N_3
 }{N_2N_1^{\frac{1}{2}}N_2^{\frac{1}{2}}N_2\langle N_3\rangle\langle
 N_4\rangle \langle N_5\rangle^{\frac{1}{2}-}\cdots \langle N_{2k+2}\rangle^{\frac{1}{2}-}}\|I\phi_1\|_{X^{\delta}_{0,\frac{1}{2}+}}\prod\limits^{2k+2}_{i=2}\|I\phi_i\|_{X^{\delta}_{1,\frac{1}{2}+}}\\
 &\lesssim& N^{-2+}
 N^{0-}_{max}\|I\phi_1\|_{X^{\delta}_{0,\frac{1}{2}+}}\prod\limits^{2k+2}_{i=2}\|I\phi_i\|_{X^{\delta}_{1,\frac{1}{2}+}},
 \end{eqnarray*}
where in the second inequality we have used Sobolev embedding \eqref{LI} to bound the terms with $j\geq 5$.\\

{\bf Case 3:} $N_2\gg N_3 \gtrsim N$.\\

We use in this instance a
trivial pointwise bound on the symbol
\begin{equation}\label{MULT2}
\left|1-\frac{m(\xi_2+\xi_3+\cdots+\xi_{2k+2})}{m(\xi_2)m(\xi_3)\cdots
m(\xi_{2k+2})}\right|\lesssim\frac{m(\xi_1)}{m(\xi_2)m(\xi_3)\cdots
m(\xi_{2k+2})}
\end{equation}
Since $m(N_1)\thicksim m(N_2)$, applying Lemma \ref{L3} we have
\begin{eqnarray*}
 {\bf Term} &\lesssim& \frac{N_1}{m(N_3)m(N_4)\cdots m(N_{2k+2})}\|I\phi_1I\phi_3\|_{L^2(\mathbb{R}^2\times
 [0,\delta])}\|I\phi_2I\phi_4\|_{L^2(\mathbb{R}^2\times[0,\delta])}\prod\limits_{j=5}^{2k+2}\|I\phi_j\|_{L^{\infty}(\mathbb{R}\times[0,\delta])}\\
 &\lesssim&\frac{N_1
 }{m(N_3)m(N_4)\cdots m(N_{2k+2})N_1^{\frac{1}{2}}N_2^{\frac{1}{2}}N_2N_3\langle
 N_4\rangle\langle N_5\rangle^{\frac{1}{2}-}\cdots \langle N_{2k+2}\rangle^{\frac{1}{2}-}}\|I\phi_1\|_{X^{\delta}_{0,\frac{1}{2}+}}\prod\limits^{2k+2}_{i=2}\|I\phi_i\|_{X^{\delta}_{1,\frac{1}{2}+}}\\
&\lesssim&\frac{1
 }{m(N_3)N_3m(N_4)\langle
 N_4\rangle m(N_{5})\langle N_{5}\rangle^{\frac{1}{2}-} \cdots m(N_{2k+2})\langle N_{2k+2}\rangle^{\frac{1}{2}-} N_2}\|I\phi_1\|_{X^{\delta}_{0,\frac{1}{2}+}}\prod\limits^{2k+2}_{i=2}\|I\phi_i\|_{X^{\delta}_{1,\frac{1}{2}+}}\\
 &\lesssim& N^{-2+}
 N^{0-}_{max}\|I\phi_1\|_{X^{\delta}_{0,\frac{1}{2}+}}\prod\limits^{2k+2}_{i=2}\|I\phi_i\|_{X^{\delta}_{1,\frac{1}{2}+}}.
 \end{eqnarray*}
where in the last inequality we use the fact that, for any
$p>0$, such that $s+p\geq 1$ the function $m(x)x^p$ is increasing and
$m(x)\langle x\rangle^p$ is bounded below, which implies
$m(N_3)N_3\gtrsim m(N)N=N$ and $m(N_4)\langle N_4\rangle\gtrsim1$,
$m(N_i)\langle N_i\rangle^{\frac{1}{2}-}\gtrsim1$, $i=5,\cdots,2k+2$.\\

{\bf Case 4:} $N_2\sim N_3\gtrsim N$ and $N_3\geq N_4$.\\

The condition $\sum_{i=1}^{2k+2}\xi_i=0$ implies $N_1\lesssim N_2$. We again bound the multiplier (\ref{MULT}) pointwise by (\ref{MULT2}). 
To obtain the decay $N^{-2+}$ we split this case into four subcases.\\

{\bf Case 4.(a):} $N_4\gtrsim N$ and $N_4\ll N_3$.\\

From (\ref{MULT2}), \eqref{L6} and Lemma \ref{L3}, we have that
\begin{eqnarray*}
 {\bf Term} &\lesssim& \frac{N_1m(N_1)}{m(N_2)m(N_3)\cdots m(N_{2k+2})}\prod\limits_{i=\{1,3\}}\|I\phi_i\|
 _{L^4(\mathbb{R}^2\times[0,\delta])}\|I\phi_2I\phi_4\|_{L^2(\mathbb{R}^2\times[0,\delta])}\prod\limits_
 {j=5}^{2k+2}\|I\phi_j\|_{L^{\infty}(\mathbb{R}\times[0,\delta])}\\
 &\lesssim&\frac{N_1m(N_1)
 }{m(N_2)m(N_3)\cdots m(N_{2k+2})N_2^{\frac{1}{2}}N_2N_3N_4\langle N_5\rangle^{\frac{1}{2}-}\cdots \langle N_{2k+2}\rangle^{\frac{1}{2}-}}\|I\phi_1\|_{X^{\delta}_{0,\frac{1}{2}+}}\prod\limits^{2k+2}_{i=2}\|I\phi_i\|_{X^{\delta}_{1,\frac{1}{2}+}}\\
&\lesssim&\frac{N_{max}^{0-}
 }{m(N_2)N_2^{\frac{3}{4}-}m(N_3)N_3^{\frac{3}{4}}m(N_4)N_4}\|I\phi_1\|_{X^{\delta}_{0,\frac{1}{2}+}}\prod\limits^{2k+2}_{i=2}\|I\phi_i\|_{X^{\delta}_{1,\frac{1}{2}+}}\\
 &\lesssim& N^{-\frac{5}{2}+}
 N^{0-}_{max}\|I\phi_1\|_{X^{\delta}_{0,\frac{1}{2}+}}\prod\limits^{2k+2}_{i=2}\|I\phi_i\|_{X^{\delta}_{1,\frac{1}{2}+}}\\
 &\lesssim& N^{-2+}
 N^{0-}_{max}\|I\phi_1\|_{X^{\delta}_{0,\frac{1}{2}+}}\prod\limits^{2k+2}_{i=2}\|I\phi_i\|_{X^{\delta}_{1,\frac{1}{2}+}}.
 \end{eqnarray*}

{\bf Case 4.(b):} $N_4\gtrsim N$ and $N_4\sim N_3$.\\

Applying the same arguments as above
\begin{eqnarray*}
 {\bf Term} &\lesssim& \frac{N_1m(N_1)}{m(N_2)m(N_3)\cdots m(N_{2k+2})}\prod\limits_{i=1}^{4}\|I\phi_i\|
 _{L^4(\mathbb{R}^2\times[0,\delta])}\prod\limits_
 {j=5}^{2k+2}\|I\phi_j\|_{L^{\infty}(\mathbb{R}\times[0,\delta])}\\
 &\lesssim&\frac{N_1m(N_1)
 }{m(N_2)m(N_3)\cdots m(N_{2k+2})N_2N_3N_4\langle N_5\rangle^{\frac{1}{2}-}\cdots \langle N_{2k+2}\rangle^{\frac{1}{2}-}}\|I\phi_1\|_{X^{\delta}_{0,\frac{1}{2}+}}\prod\limits^{2k+2}_{i=2}\|I\phi_i\|_{X^{\delta}_{1,\frac{1}{2}+}}\\
&\lesssim&\frac{N_{max}^{0-}
 }{m(N_2)N_2^{\frac{2}{3}-}m(N_3)N_3^{\frac{2}{3}}m(N_4)N_4^{\frac{2}{3}}}\|I\phi_1\|_{X^{\delta}_{0,\frac{1}{2}+}}\prod\limits^{2k+2}_{i=2}\|I\phi_i\|_{X^{\delta}_{1,\frac{1}{2}+}}\\
  &\lesssim& N^{-2+}
 N^{0-}_{max}\|I\phi_1\|_{X^{\delta}_{0,\frac{1}{2}+}}\prod\limits^{2k+2}_{i=2}\|I\phi_i\|_{X^{\delta}_{1,\frac{1}{2}+}}.
 \end{eqnarray*}

{\bf Case 4.(c):} $N_4\ll N$ and $N_1\ll N_2$.\\

Again using the bound (\ref{MULT2}) and Lemma \ref{L3}, we have
\begin{eqnarray*}
{\mathbf{Term}} \!\!\!&\lesssim& \!\!\! \frac{N_1m(N_1)}{m(N_2)\cdots m(N_{2k+2})} 
\left\|{I\phi_1}I\phi_2\right\|_{L^2(\R\times[0,\delta])} \left\|{I\phi_3}I\phi_4\right\|_{L^2(\R\times[0,\delta])} 
\prod\limits_{j=5}^{2k+2}\|I\phi_j\|_{L^{\infty}(\mathbb{R}\times[0,\delta])}\\
&\lesssim& \!\!\! \frac{N_1m(N_1)} 
{m(N_2)m(N_3)m(N_4)N_2^{\frac{1}{2}}N_3^{\frac{1}{2}} N_2N_3 \langle N_4\rangle 
\langle N_5\rangle^{\frac{1}{2}-}\cdots \langle N_{2k+2}\rangle^{\frac{1}{2}-}}
\|I\phi_1\|_{X^{\delta}_{0,\frac{1}{2}+}} \prod_{i=2}^{2k+2} \|I\phi_i\|_{X^{\delta}_{1,\frac{1}{2}+}}\\
&\lesssim&  \!\!\!\frac{N_{max}^{0-}} {m(N_2)N_2^{1-}m(N_3)N_3m(N_4) \langle N_4\rangle}
\|I\phi_1\|_{X^{\delta}_{0,\frac{1}{2}+}} \prod_{i=2}^{2k+2} \|I\phi_i\|_{X^{\delta}_{1,\frac{1}{2}+}}\\
&\lesssim&  \!\!\!N^{-2+}N_{max}^{0-}\|I\phi_1\|_{X^{\delta}_{0,\frac{1}{2}+}} \prod_{i=2}^{2k+2} \|I\phi_i\|_{X^{\delta}_{1,\frac{1}{2}+}}.
\end{eqnarray*}

{\bf Case 4(d):} $N_4\ll N$ and $N_1\sim N_2\sim N_3 \gtrsim N$.\\

In this case, we use an argument inspired by Pecher \cite[Proposition 5.1]{P1}. Since $\sum_{i=1}^{2k+2}\xi_i=0$, two of the large frequencies have different sign,
say, $\xi_1$ and $\xi_2$. Indeed, if all have the same size, we obtain $|\xi_1+\xi_2+\xi_3|\geq |\xi_3| \gg |\xi_4+\cdots+\xi_{2k+2}|$, a contradiction with $\sum_{i=1}^{2k+2}\xi_i=0$). Thus,
\begin{equation*}
|\xi_1|\leq |\xi_1-\xi_2|\leq 2|\xi_1|
\end{equation*}
and
\begin{equation*}
|\xi_1+\xi_2|= |\xi_3+\xi_4|\sim|\xi_1|.
\end{equation*}

Therefore, using the bound \eqref{MULT2} and Lemma \ref{L3}, we have
\begin{eqnarray*}
{\mathbf{Term}} \!\!\!&\lesssim& \!\!\!
\frac{N_1^{\frac{1}{2}}m(N_1)}{m(N_2)\cdots m(N_{2k+2})} 
\left\|(D_x^{\frac{1}{2}}I\phi_1)I\phi_2\right\|_{L^2_{x,t}} \left\|I\phi_3I\phi_4\right\|_{L^2_{x,t}} \prod\limits_
 {j=5}^{2k+2}\|I\phi_j\|_{L^{\infty}(\mathbb{R}\times[0,\delta])}\\
&\lesssim& \!\!\! \frac{N_1^{\frac{1}{2}}m(N_1)} 
{m(N_2)m(N_3)m(N_4)N_3^{\frac{1}{2}} N_3 \langle N_4\rangle\langle N_5\rangle^{\frac{1}{2}-}\cdots \langle N_{2k+2}\rangle^{\frac{1}{2}-}}
\left\|(D_x^{\frac{1}{2}}I\phi_1)I\phi_2\right\|_{L^2_{x,t}} 
\prod_{i=3}^{4} \|I\phi_i\|_{X^{\delta}_{1,\frac{1}{2}+}} \prod\limits_
 {j=5}^{2k+2}\|I\phi_j\|_{X^{\delta}_{1,\frac{1}{2}+}}\\
&\lesssim& \!\!\! \frac{N_{max}^{0-}} {m(N_2)N_2^{1-}m(N_3)N_3m(N_4) \langle N_4\rangle}\|I\phi_1\|_{X^{\delta}_{0,\frac{1}{2}+}} 
\prod_{i=2}^{2K+2} \|I\phi_i\|_{X^{\delta}_{1,\frac{1}{2}+}}\\
&\lesssim& \!\!\! N^{-2+}N_{max}^{0-}\|I\phi_1\|_{X^{\delta}_{0,\frac{1}{2}+}} \prod_{i=2}^{2k+2} \|I\phi_i\|_{X^{\delta}_{1,\frac{1}{2}+}},
\end{eqnarray*}
where we have estimated $\|(D_x^{\frac{1}{2}}I\phi_1)I\phi_2\|_{L^2_{x,t}}$ via Lemma \ref{L3}.\\
\fim

\section{Global theory}\label{GWP}

Before proceeding to the proof of Theorem \ref{T1} we need to establish a variant of local well-posedness result for the modified equation \eqref{MNLB}.
%\begin{eqnarray}\label{MBOU}
%\left\{
%\begin{array}{c}
%Iu_{tt}-I(u_{xx})+I(u_{xxxx})=I((|u|^{2k}u)_{xx})\\
%Iu(x,0)=I\phi(x);\ \ \partial_tIu(x,0)=I\psi(x)_x.
%\end{array}\right.
%\end{eqnarray}
Clearly if $Iu\in H^1(\R)$ is a solution of \eqref{MNLB}, then $u\in H^s(\R)$ is a solution of \eqref{NLB} in the same time interval. 

Next we prove a local existence result for this modified equation. Since we do not have
scaling invariance we also need to estimate the solution existence time. The crucial nonlinear estimate for the local existence 
is given in the next lemma.
\begin{lemm}\label{NLEL}
If $s>\dfrac{4k-5}{4(2k-1)}$, $k\in\mathbb{N}$, then
\begin{equation}\label{NLE}
\||u|^{2k}u\|_{X_{s,0}}\lesssim
\|u\|_{X_{s,\frac{1}{2}+}}^{2k+1},
\end{equation}
\end{lemm}

\noindent\textbf{Proof. } It is easy to see that, for all $s>0$
$$\langle \xi_1+\cdots+\xi_{2k+1}\rangle^s\lesssim \langle \xi_1\rangle^s +\cdots+\langle \xi_{2k+1}\rangle^s \peq \textrm{for all} \peq s>0.$$

Thus, by duality and a Leibniz rule, \eqref{NLE} follows from
\begin{equation}
\left|\int_{\mathbb{R}}\int_{\mathbb{R}}J^s\psi_1\prod\limits_{i=2}^{2k+2}\psi_i
dxdt\right|\lesssim
\left(\prod\limits_{i=1}^{2k+1}\|\psi_i\|_{X_{s,\frac{1}{2}+}}\right)\|\psi_{2k+2}\|_{X_{0,0}}.
\end{equation}
First, we use H\"{o}lder's inequality on the left hand side of \eqref{NLE}, taking the
factors in
$L_{x,t}^6,L_{x,t}^{12},L_{x,t}^{4(2k-1)},\cdots,L_{x,t}^2$. Thus, applying the Sobolev embedding and the Strichartz inequality \eqref{L81}, we have
\begin{eqnarray*}
\|J^s\psi_1\|_{L_{x,t}^6(\mathbb{R}^{1+1})}&\lesssim & \|J^s\psi_1\|_{X_{0,\frac{1}{2}+}}\\
&= &\|\psi_1\|_{X_{s,\frac{1}{2}+}},
\end{eqnarray*}
\begin{eqnarray*}
\|\psi_2\|_{L_{x,t}^{12}(\mathbb{R}^{1+1})}&\lesssim &
\|J^{\frac{1}{4}}\psi_2\|
_{L_t^{12}L_{x}^3(\mathbb{R}^{1+1}))}\\
&\lesssim &\|J^{\frac{1}{4}}\psi_2\|_{X_{0,\frac{1}{2}+}}\\
 &\lesssim &\|\psi_2\|_{X_{s,\frac{1}{2}+}},
\end{eqnarray*}
\begin{eqnarray*}
\|\psi_3\|_{L_{x,t}^{4(2k-1)}(\mathbb{R}^{1+1})}&\lesssim &
\|J^{\frac{4k-5}{4(2k-1)}}\psi_3\|
_{L_t^{4(2k-1)}L_{x}^{\frac{2k-1}{k-1}}(\mathbb{R}^{1+1}))}\\
&\lesssim &\|J^{\frac{4k-5}{4(2k-1)}}\psi_3\|_{X_{0,\frac{1}{2}+}}\\
 &\lesssim &\|\psi_3\|_{X_{s,\frac{1}{2}+}}.
\end{eqnarray*}

By similar arguments as the previous one above, for $j=4,5,\cdots,2k+1$, we obtain
\begin{eqnarray*}
\|\psi_j\|_{L_{x,t}^{4(2k-1)}(\mathbb{R}^{1+1})}
 \lesssim \|\psi_j\|_{X_{s,\frac{1}{2}+}}.
\end{eqnarray*}

Finally, applying the trivial estimate \eqref{TE} we have
$$\|\psi_{2k+2}\|_{L^2_tL^2_x}\leqslant
\|\psi_{2k+2}\|_{X_{0,0}}.$$ 

Therefore, the above inequalities together with the fact that $\dfrac{1}{4}<\dfrac{4k-5}{2(2k-1)}<\dfrac{1}{2}$ for all $k>2$ yield \eqref{NLE}.\\
\fim
 
\begin{rema}
It should be interesting to prove inequality \eqref{NLE} for $s>{1}/{2}-{1}/{k}$. As a consequence, one can recover all the well known range of existence 
for the local theory given in Farah \cite{LG3} in terms of the $X_{s,b}$ spaces.
\end{rema}

Applying the interpolation lemma (see \cite{CKSTT2}, Lemma 12.1) to \eqref{NLE} we obtain
\begin{equation*}
\|I(|u|^{2k}u)\|_{X_{1,0}}\lesssim
\|u\|_{X_{1,\frac{1}{2}+}}^{2k+1}.
\end{equation*}
where the implicit constant is independent of $N$. Now standard arguments invoking the contraction-mapping principle give the 
following variant of local well posedness.

\begin{theo}\label{t3.2}

Assume $s<1$, $(\phi,\psi)\in H^s(\mathbb{R})\times
H^{s-1}(\mathbb{R})$ be given. Then there exists a positive number
$\delta$ such that IVP (\ref{MNLB}) has a unique local solution $Iu\in
C([0,\delta], H^1(\mathbb{R}))$ such that
\begin{equation}
\max\{\|(-\Delta)^{-\frac{1}{2}}\partial_tIu\|_{X^{\delta}_{0,\frac{1}{2}+}}, \|Iu\|_{X_{1,\frac{1}{2}+}^{\delta}}\}\leqslant
C(\|I\phi\|_{H^1}+\|I\psi\|_{L^2}).
\end{equation}
Moreover, the existence time can be estimates by
\begin{equation}\label{MDE}
\delta^{\frac{1}{2}-}\thicksim \frac{1}{(\|I\phi\|_{H^1}+\|I\psi\|_{L^2})^{2k}}.
\end{equation}
\end{theo}

We should note that the power $\frac{1}{2}-$ in \eqref{MDE} is also closely related to the index $s$ obtained in our global result. Since we need to iterate at the very end of the method, we would like to maximize the time of existence $\delta$ (given in the above theorem) in each step of iteration. Since the denominator in the right hand side of  \eqref{MDE} is very large the only way to do that is to maximize the power of $\delta$ in the left hand side. By standards linear estimates in $X_{s,b}$ spaces (see \cite[Lemma 4.2]{LGF}) the best that we can obtain is $\frac{1}{2}-$.

Now, we have all tools to prove our global result stated in Theorem \ref{T1}.\\

\noindent\textbf{Proof of Theorem \ref{T1}. } Let $(\phi, \psi) \in H^s(\R) \times H^{s-1}(\R)$
with $1/2\leq s <1$. 
Our goal is to construct a solution to \eqref{MNLB} (and therefore to \eqref{NLB}) on an arbitrary time interval $[0,T]$. 
From the definition of the multiplier $I$ we have
\begin{eqnarray*}
\|I\phi\|_{H^1}&\leq& c N^{1-s}\|\phi\|_{H^{s}},\\
\|I\psi\|_{L^2}&\leq& c N^{1-s}\|\psi\|_{H^{s-1}}.
\end{eqnarray*}

Therefore, there exists $c_1>0$ such that
\begin{equation*}
\max\left\{\|I\phi\|_{H^1}, \|I\psi\|_{L^2}\right\}\leq c_1 N^{1-s}.
\end{equation*}

We use our local existence theorem on $[0, \delta]$, where $\delta^{\frac{1}{2}-}\sim N^{-2k(1-s)}$ and conclude
\begin{eqnarray}\label{IU0}
\max\left\{\|Iu\|_{X^{\delta}_{1,\frac{1}{2}+}}, \|(-\Delta)^{-\frac{1}{2}}\partial_tIu\|_{X^{\delta}_{0,\frac{1}{2}+}} \right\}&\leq& c\left(\|I\phi\|_{H^1}+\|I\psi\|_{L^2}\right)\\
&\leq& c_2 N^{1-s}.
\end{eqnarray}

From the conservation law (\ref{EC}), we obtain
\begin{equation}\label{CC2}
\|Iu(\delta)\|^2_{H^1}+\|(-\Delta)^{-\frac{1}{2}}\partial_tIu(\delta)\|^2_{L^2}\leq c_3 E(Iu(\delta)).
\end{equation}

On the other hand, since $\|f\|_{L^{2k+2}}\lesssim \|f\|_{H^s}$, for $s>\dfrac{k}{2(k+1)}$ (note that $\dfrac{k}{2(k+1)}<\dfrac{1}{2}$), 
we have
\begin{equation*}
E(Iu(0))\leq c \,N^{2(1-s)}+c\|\phi\|_{L^{2k+2}}^{2k+2}\leq c_4N^{2(1-s)}.
\end{equation*} 

By the almost conservation law stated in Proposition \ref{p4.1} and \eqref{IU0}, we have
\begin{equation*}
E(Iu(\delta))\leq E(Iu(0))+cN^{-2+}N^{4(1-s)}<2c_4N^{2(1-s)}.
\end{equation*}

We iterate this process $M$ times obtaining 
\begin{equation}\label{UB}
E(Iu(\delta))\leq E(Iu(0))+cMN^{-2+}N^{4(1-s)}<2c_4N^{2(1-s)}.,
\end{equation}
as long as $cMN^{-2+}N^{4(1-s)}<c_4N^{2(1-s)}$, which implies that the lifetime of the local results remains uniformly of size 
$\delta^{\frac{1}{2}-}\sim N^{-2k(1-s)}$. 

Given a time $T>0$, the number of iteration steps to reach this time is $T\delta^{-1}$. Therefore, to carry out $T\delta^{-1}$
iterations on time intervals, before the quantity $E(Iu)(t)$ doubles, the following condition has to be
fulfilled
\begin{eqnarray}\label{GGV2}
N^{-2+}N^{(2k+2)(1-s)}T\delta^{-1}\ll N^{2-2s}.
\end{eqnarray}

Since $\delta^{\frac{1}{2}-}\sim N^{-2k(1-s)}$, the condition (\ref{GGV2}) can be obtained for 
\begin{equation}\label{expot}
T\sim N^{(6ks-6k+2)-}.
\end{equation}

\begin{rema}
Note that the exponent of $N$ on the right hand side of (\ref{expot}) is positive provided $s>1-({1}/{3k})$, 
hence the definition of $N$ makes sense for arbitrary large $T$.
\end{rema}

Finally, we need to establish the polynomial bound (\ref{pb}). By our choice of $N$, relation (\ref{smo}) and (\ref{CC2}) 
imply for $T\gg 1$ that
\begin{equation*}
\sup_{t\in[0,T]}\left\{\|u(t)\|^2_{H^{s}}, \|(-\Delta)^{-\frac{1}{2}}\partial_tu(t)\|^2_{H^{s-1}}\right\}\lesssim E(Iu(T)) 
\lesssim N^{2(1-s)} \sim T^{\frac{2(1-s)}{6ks-6k+2}+}
\end{equation*}
which implies the polynomial bound (\ref{pb}).\\
\fim

%%%%%%%%%%%%%%%%%%%%%%%%%%%%%%%%%%%%%%%%%%%%%%%%%%%%%%%%%%%%%%%%%%%%%%%%%%%%%%
\centerline{\textbf{Acknowledgment}}

We thank Felipe Linares for detailed comments and corrections.
%%%%%%%%%%%%%%%%%%%%%%%%%%%%%%%%%%%%%%%%%%%%%%%%%%%%%%%%%%%%%%%%%%%%%%%%%%%%%%
%
% \bibliographystyle{abbrv}
% \bibliography{bibfarah}
% \bibliographystyle{abbrv}

\end{document}